\documentstyle[a4,11pt]{article}
\begin{document}

\begin{titlepage}
\vskip 2cm
\begin{flushright}
Preprint CNLP-1994-06
\end{flushright}
\vskip 2cm
\begin{center}
{\bf INTEGRABLE GEOMETRY and  SOLITON EQUATIONS IN 2+1 DIMENSIONS}\footnote{Preprint
CNLP-1994-06. Alma-Ata.1994 }
\end{center}
\vskip 2cm
\begin{center}
{\bf R.Myrzakulov}
\end{center}

\vskip 1cm
Centre for Nonlinear Problems, PO Box 30, 480035, Alma-Ata-35, Kazakhstan\\
E-mail: cnlpmyra@satsun.sci.kz

\vskip 1cm

\begin{abstract}

Using the differential geometry of curves and surfaces,  the L-equivalent
soliton equations of the some (2+1) - dimensional integrable spin systems
are found. These equations include the modified Novikov-Veselov, Kadomtsev-
Petviashvili, Nizhnik-Novikov-Veselov and other equations.
Some aspects of the connection between geometry and
multidimensional soliton equations are discussed.

\end{abstract}


\end{titlepage}

\setcounter{page}{1}
\newpage
\vskip 2cm

\begin{center}
\section{Introduction}
\end{center}

The relation between geometry and soliton equations
has been the subject of continued interest in recent years[1-9]. In general
the connection between geometry and nonlinear partial differential
equations(NPDE) has the  long history (for a historical review, see, e.g.,
[4]). Spin systems (which are the important from physical and mathematical
point of views subclass of NPDE) are a good laboratory to demonstrate and to
understand the connection of geometry and NPDE. The first representative of
integrable spin systems is the isotropic Landau-Lifshitz equation(LLE)[2,10]
$$
\vec S_{t} = \vec S \wedge \vec S_{xx}   \eqno(1)
$$
where $\vec S^{2} = E = \pm 1$. Here and hereafter subcripts denote
partial derivatives. Pioneering work by Lakshmanan [2] showed that the LLE
(1) is equivalent to the known nonlinear Schrodinger equation(NLSE)
$$
iq_t + q_{xx} + 2E\mid q \mid ^{2} q  = 0 \eqno (2)
$$
for the case $E = + 1$ (for the case $E =  -1$ the such
equivalence was established in [17]). This equivalence in [17] we called
{\bf the Lakshmanan equivalence} or shortly {\bf L-equivalence}. As well known
between equations(1) and (2) take place the gauge equivalence[15].

Let us, first, we briefly recall some of the well known facts on the
geometrical formalism that used in [2], with the minor modifications
(including the case $E = -1$) of ref.[17]. Consider the motion of curves
which are given by
$$
\left ( \begin{array}{ccc}
\vec e_{1} \\
\vec e_{2} \\
\vec e_{3}
\end{array} \right)_{x} = C
\left ( \begin{array}{ccc}
\vec e_{1} \\
\vec e_{2} \\
\vec e_{3}
\end{array} \right), \eqno(3a)
$$
$$
\left ( \begin{array}{ccc}
\vec e_{1} \\
\vec e_{2} \\
\vec e_{3}
\end{array} \right)_{t} = G
\left ( \begin{array}{ccc}
\vec e_{1} \\
\vec e_{2} \\
\vec e_{3}
\end{array} \right) \eqno(3b)
$$
with
$$
C =
\left ( \begin{array}{ccc}
0   & k     & 0 \\
-Ek & 0     & \tau  \\
0   & -\tau & 0
\end{array} \right) ,\quad
G =
\left ( \begin{array}{ccc}
0       & \omega_{3}  & -\omega_{2} \\
-E\omega_{3} & 0      & \omega_{1} \\
E\omega_{2}  & -\omega_{1} & 0
\end{array} \right).
$$
Hence we have
$$
C_t - G_x + [C, G] = 0
  \eqno (4)
$$
Here $\vec e_{1}, \vec e_{2}, \vec e_{3}$ denote respectively
the unit tangent, normal and binornal vectors, defined in the usual
way, $k$ and $\tau$ are respectively the curvature and torsion
of the curve. Note that (3a) is the usual Serret-Frenet equation
(SFE). If $\vec S = \vec e_{1}$, then the
function $ q = \frac{k}{2}e^{-i\partial^{-1}_{x} \tau}$ satisfies the
NLSE(2), that is, equations(1) and (2) are L-equivalent each to other.

As known recently many efforts have been made to study the (2+1)-dimensional
integrable NPDE[11-14].
Here we have the interesting fact: the (1+1) - dimensional integrable
NPDE admit some number (not one) integrable (and nonintegrable) (2+1) - dimensional
generalizations. So, for example, the LLE(1) has the following
(2+1) - dimensional integrable and nonintegrable extensions: \\
$1^{\circ}$. The Myrzakulov I (M-I) equation[17]
$$ \vec S_{t}=(\vec S\wedge \vec S_{y}+u\vec S)_x \eqno (5a) $$
$$ u_x=-\vec S_(\vec S_{x}\wedge \vec S_{y}) \eqno (5b) $$
$2^{\circ}$. The Myrzakulov VIII (M-VIII) equation[17]
$$ iS_t=\frac{1}{2}[S_{xx},S]+iwS_{x}   \eqno (6a) $$
$$ w_{y}=\frac{1}{4i}tr(S[S_{x},S_{y}])  \eqno (6b) $$
$3^{\circ}$. The Ishimori   equation[17]
$$ iS_t+\frac{1}{2}[S,(\frac{1}{4}S_{xx}+\alpha^{2}S_{yy})]+
iu_{y}S_x+iu_{x}S_y = 0 \eqno(7a)$$
$$ \alpha^{2}u_{yy} - \frac{1}{4}u_{xx}=
\frac{\alpha^{2}}{4i}tr(S[S_y,S_x]) \eqno(7b)$$
$4^{\circ}$. The Myrzakulov IX (M-IX) equation[17]
$$ iS_t+\frac{1}{2}[S,M_1S]+A_2S_x+A_1S_y = 0 \eqno(8a)$$
$$ M_2u=\frac{\alpha^{2}}{4i}tr(S[S_y,S_x]) \eqno(8b)$$
$5^{\circ}$. The Myrzakulov XVIII (M-XVIII) equation[17]
$$ iS_t+\frac{1}{2}[S,(\frac{1}{4}S_{xx}-\alpha(2b+1)S_{xy}+\alpha^{2}S_{yy})]+
A_{20}S_x+A_{10}S_y = 0 \eqno(9a)$$
$$ \alpha^{2}u_{yy} - \frac{1}{4}u_{xx}=\frac{\alpha^{2}}{4i}tr(S[S_y,S_x]) \eqno(9b)$$
$6^{\circ}$. The (2+1) - dimensional LLE
$$
\vec S_{t} = \vec S \wedge (\vec S_{xx} + \vec S_{yy})   \eqno(10)
$$
All of these equations in 1+1 dimension reduce to
the LLE(1). Note that here the
Ishimori(7), M-I(5), M-VIII(6), M-IX(8) and M-XVIII(9) equations are integrable,
at the same time equation (10) is not integrable.

In the study of the (2+1) - dimensional NPDE naturally arise the
following questions:

$1^{\circ}$. {\it What is the analog of the SFE(3a) in 2+1 dimensions?}

{\bf Answer:  The analog of the SFE(3a) in 2+1 dimensions is the following
set of equations:}
$$
\left ( \begin{array}{ccc}
\vec e_{1} \\
\vec e_{2} \\
\vec e_{3}
\end{array} \right)_{x} = C
\left ( \begin{array}{ccc}
\vec e_{1} \\
\vec e_{2} \\
\vec e_{3}
\end{array} \right), \quad
\left ( \begin{array}{ccc}
\vec e_{1} \\
\vec e_{2} \\
\vec e_{3}
\end{array} \right)_{y} = D
\left ( \begin{array}{ccc}
\vec e_{1} \\
\vec e_{2} \\
\vec e_{3}
\end{array} \right) \eqno(11)
$$
with
$$
C =
\left ( \begin{array}{ccc}
0   & k     & 0 \\
-Ek & 0     & \tau  \\
0   & -\tau & 0
\end{array} \right) ,\quad
D =
\left ( \begin{array}{ccc}
0        & m_{3}  & m_{2} \\
-Em_{3}  & 0      & m_{1} \\
E m_{2}  & -m_{1} & 0
\end{array} \right). \eqno(12)
$$
In [17] , this set of equations was called the
Serret-Frenet-Myrzakulov equation(SFME). [ Compare the SFME(11) with the
usual SFE(3a)].
Note that from the SFME(11)
follows the following conditions
$$
C_y - D_x + [C, D] = 0,  \eqno (13a)
$$
or in terms of elements
$$
m_{1} = \partial ^{-1}_{x}(\tau_{y} + Ek m_{2}), \eqno(13b)
$$
$$
m_{2} = \partial ^{-1}_{x}(k m_{1} - \tau m_{3}) \eqno(13c)
$$
$$
m_{3}  = \partial ^{-1}_{x}(k_{y} - \tau  m_{2}).  \eqno(13d)
$$
These conditions play a key role in the geometrical interpretation
of the (2+1) - dimensional NPDE from the integrability point of view.
In particular, as it seems to us, these conditions suggest to us,
how find the integrable reductions of the (2+1) - dimensional NPDE.
Let us return to questions.

$2^{\circ}$. {\it What is the (2+1)-dimensional
version of the  L-equivalence?}

$3^{\circ}$. {\it How find the L-equivalent counterparts of the known and new
integrable and nonintegrable NPDE?}\\
and so on. May be the some answers of the some these questions were given in [17].
It should be noted
that the  SFME(11) plays a key role in our construction.
The other new moment of our formalism, that is, it contain both cases: the
focusing version when $E = +1$ (the Euclidean case) and the defocusing case
when $E = -1$ (the Minkowski case).

Third, the central element of our construction is
the Myrzakulov - 0(M-0) equation. Let us consider a $n$ - dimensional space
$R^{n}$ with the unit basic vectors $e_{j},  j = 1,2, ... ,n.$
Let $\vec S \equiv  \vec e_{1}, \vec S^{2} \equiv \vec e_{1}^{2} = E = \pm 1$.
Then, for example, the (3+1) - dimensional M-0 equation reads as[9]
$$
\vec S_{t} = \sum_{j=2}^{n} a_{j} \vec e _{j}   \eqno(14a)
$$
$$
\vec S_{x} = \sum_{j=2}^{n} b_{j} \vec e _{j}   \eqno(14b)
$$
$$
\vec S_{y} = \sum_{j=2}^{n} c_{j} \vec e _{j}   \eqno(14c)
$$
$$
\vec S_{z} = \sum_{j=2}^{n} d_{j} \vec e _{j}   \eqno(14d)
$$

In this paper we consider the some integrable reductions of
the (2+1) - dimensional M-0 equation. Using the geometrical formalism
the L-equivalent
counterparts of the some (2+1) - dimensional integrable spin systems
were found.

\section{The Myrzakulov XVII equation}

The Myrzakulov XVII(M-XVII) equation
$$
\vec S_t = \frac{1}{4}\vec S_{xxx} - \frac{3}{4}\vec S_{xyy} +
C_{1}\vec S_{x} + C_{2}\vec S_{y}
+ C_{3}\vec S,  \eqno (15a)
$$
$$
V_{x} + iV_{y} = \frac{1}{8}[(
\vec S_{x}^{2} + \vec S_{y}^2)_{x} -
i(\vec S_{x}^{2} + \vec S_{y}^2)_{y}]
\eqno (15b)
$$
was introduced in [17] and arises from the compatibility conditions of
the linear problem
$$
\Phi_{y} = S\Phi_{x},\quad
\Phi_{t} = \Phi_{xxx} + B_{2}\Phi_{xx} + B_{1}\Phi_{x}  \eqno(16)
$$
Here
$$
B_{2} = \frac{3}{2}SS_{x}, \quad
B_{1} = \frac{3}{4}SS_{xx} + \frac{3}{16}(\vec S_{x}^{2} + \vec S_{y}^{2}
+8\bar V + 8V)I + \frac{3}{2}(\bar V - V)S + \frac{3i}{8}S\{S_{x},S_{y}\}
+ \frac{3i}{4}S_{xy}
$$
$$
C_{1} = \frac{3}{16}(\vec S_{x}^{2} + \vec S_{y}^{2} +8\bar V + 8V),\,\,\,\,\,
C_{2} =  \frac{3}{4}(2i\bar V - 2iV - \vec S_{x}\vec S_{y}),
$$
$$
C_{3} = - \frac{3}{2}\{[\bar V + V - \frac{3}{8}\vec S_{x}^{2}]_{x} + \frac{1}{2}
(\vec S_{x} \vec S_{y})_{y}\}S = \frac{1}{4}\vec S (3\vec S_{xyy} - \vec S_{xxx}),
$$
$\vec S =  (S_{1}, S_{2})$ is the spin vector,
$ \vec S^{2} = E = \pm 1,\quad V$ is scalar function, \quad
$
S = \sum _{k=1}^{3}S_{k}\sigma _{k}(S_{3} = 0),\quad \sigma_{k}$
are Pauli matrix,
$$
\sigma_{1} =
\left( \begin{array}{cc}
0            &  1  \\
1            &  0
\end{array} \right),\,\,\,
\sigma_{2} =
\left( \begin{array}{cc}
0           &  -i  \\
i           &  0
\end{array} \right),\,\,\,
\sigma_{3} =
\left( \begin{array}{cc}
1           &  0  \\
0           & -1
\end{array} \right).
$$

It is convenient sometimes the following form of the M-XVII equation
$$
\vec S_t = \vec S_{zzz} + \vec S_{\bar z \bar z \bar z} + C^{-}\vec S_{z}
+ C^{+}\vec S_{\bar z}
+ C_{3}\vec S,  \eqno (17a)
$$
$$
V_{\bar z}  = \frac{1}{2}(\vec S_{z}\vec S_{\bar z})_{z}
\eqno (17b)
$$
where $ z = x + i y $.

\section{The Serret-Frenet-Myrzakulov equation
and the mNVE as the L-equivalent counterpart of the M-XVII equation}

In this section we find the L-equivalent counterpart
of the M-XVII equation(15). To this purpose, following[17] we consider the
motion of 2-dimensional curves. The Serret-Frenet-Myrzaluov equation(SFME)
in this case  has the form[17]
$$
\left ( \begin{array}{ccc}
\vec e_{1} \\
\vec e_{2}
\end{array} \right)_{x} = C
\left ( \begin{array}{ccc}
\vec e_{1} \\
\vec e_{2}
\end{array} \right),\quad
\left ( \begin{array}{ccc}
\vec e_{1} \\
\vec e_{2}
\end{array} \right)_{y} = D
\left ( \begin{array}{ccc}
\vec e_{1} \\
\vec e_{2}
\end{array} \right), \eqno(18)
$$
with
$$
C =
\left ( \begin{array}{cc}
0   & k     \\
-Ek & 0
\end{array} \right) ,\quad
D =
\left ( \begin{array}{cc}
0       & m \\
-Em     & 0
\end{array} \right), \eqno(19)
$$
where $m = \partial^{-1}_{x}k_{y}$. [ Compare the SFME(18) with the usual SFE
$$
\left ( \begin{array}{ccc}
\vec e_{1} \\
\vec e_{2}
\end{array} \right)_{x} = C
\left ( \begin{array}{ccc}
\vec e_{1} \\
\vec e_{2}
\end{array} \right), \quad
C =
\left ( \begin{array}{cc}
0   & k     \\
-Ek & 0
\end{array} \right) ]
$$
At the same time the time evolution of curve is specified by
$$
\left ( \begin{array}{cc}
\vec e_{1} \\
\vec e_{2}
\end{array} \right)_{t} = G
\left ( \begin{array}{cc}
\vec e_{1} \\
\vec e_{2}
\end{array} \right).     \eqno(20)
$$
Here $\vec e_{k}$ are a unit basic vectors,
$$
G =
\left ( \begin{array}{cc}
0        & \omega  \\
-E\omega & 0
\end{array} \right), \eqno(21)
$$
$ k$ is the curvature of curve. From these equations  we have
$$
C_y - D_x + [C, D] = 0,\quad
C_t - G_x + [C, G] = 0,\quad
D_t - G_y + [D, G] = 0. \eqno (22)
$$
So we get
$$
k_{t}=\omega_{x},\quad
m_{t}=\omega_{y}.    \eqno(23)
$$
Let $\vec S \equiv \vec e_{1}.$ Then we obtain
$$
\omega  = \frac{1}{4}(k_{xx} - 3k_{yy}) - \frac{1}{4}(k^{3} -km^{2}) +
c_{1}k + c_{2}m.   \eqno(24)
$$
with
$$
c_{1} = \frac{3}{16}(k^{2} + m^{2} + 8\bar V + 8V), \,\,\,\,\,
c_{2} = \frac{3}{4}(2i\bar V - 2iV - km).
$$
Now introduce the following new real function
$$
q = [k^{2}/4 + (\partial^{-1}_{x}k_{y})^{2}]^{\frac{1}{2}}  \eqno(25)
$$
It is not difficult check that this function satisfies the following
modified Novikov-Veselov equation(mNVE)
$$
q_t = (q_{zzz} + 3Vq_z + \frac{3}{2}V_zq) + (q_{\bar z \bar z \bar z}
+ 3 \bar V q_{\bar z} + \frac{3}{2}\bar V_{\bar z}q) \eqno (26a)
$$
$$
V_{\bar z} = (q^2)_z . \eqno (26b)
$$
where $z = x + iy$. As well known this equation was introduced in[26]
 and is associated with the Lax representation
$$
L^{mNV}\Psi =
\left( \begin{array}{cc}
\partial  &  -q  \\
q         &  \bar \partial
\end{array} \right)\Psi = 0 \eqno(27a)
$$
$$
\Psi_t = (A^{+} + A^{-})\Psi    \eqno(27b)
$$
where
$$
A^{+} = \partial ^3 + 3
\left( \begin{array}{cc}
0  &  -q_z  \\
0  &  V
\end{array} \right) \partial +
 \frac{3}{2}
\left( \begin{array}{cc}
0  &  2qV  \\
0  &  V_{z}
\end{array} \right),\quad
A^{-} = \bar \partial ^3 + 3
\left( \begin{array}{cc}
\bar V \bar \partial  &  0   \\
q_{\bar z}              &  0
\end{array} \right) \bar \partial +
 \frac{3}{2}
\left( \begin{array}{cc}
\bar V_{\bar z}  &  0   \\
-2q \bar V       &  0
\end{array} \right).
$$
So we have proved that the mNVE(26) and the M-XVII equation(15)
are L-equivalent each to other that coincide with the fact that
these equations too are gauge eqivalent each to other[16].

\section{The related spin systems and their L-equivalents}

\subsection{ The Myrzakulov XI equation}

The mNV equation(26) is the some modification of the Novikov - Veselov
equation(NVE)
$$
q_t = q_{zzz} +q_{\bar z\bar z\bar z} + (Vq)_z + (\bar Vq)_{\bar z}, \eqno (28a)
$$
$$
V_{\bar z} = 3(q^2)_{z}. \eqno (28b)
$$
The NV equation(28) was introduced in[24-25]. In contrast
with the mNVE, the NVE(28) is associated with the
$ (L, A, B) $- triple
$$ \frac{\partial}{\partial t}L^{NV} + [L^{NV},A] - BL^{NV} =0, \eqno (29) $$
где
$$ L^{NV} = \partial \bar \partial + q,\,\,\,\, A =
(\partial ^3 +V \partial) + (\bar \partial ^3 + \bar V \bar \partial ),\,\,\,\,
B= \partial V + \bar \partial \bar V. \eqno (30) $$
Let us consider the other form of the NVE
$$
q_{t} = \alpha q_{xxx} + \beta q_{yyy} - 3\alpha (vq)_{x} - 3\beta (wq)_{y} \eqno(31a)
$$
$$
w_{x} = q_{y} \eqno(31b)
$$
$$
v_{y} = q_{x}.  \eqno(31c)
$$
This equation is the compatibility condition of the linear problem
$$
\phi_{xy} = q \phi  \eqno(32a)
$$
$$
\phi_{t} = \alpha \phi_{xxx} + \beta \phi_{yyy} -
3\alpha v \phi_{x} - 3\beta w \phi_{y}  \eqno(32b)
$$
and introduced in [24-25]. Here note the NVE(31) is the  Lakshmanan equivalent to
the following (2+1)-dimensional integrable spin model - the socalled
M-XI equation[17]
$$
\vec S_{t} = \frac{\omega}{k}\vec S_{x} \eqno(33)
$$
with
$$
\omega =  \alpha q_{xxx} + \beta q_{yyy} - 3\alpha (vq)_{x} -
3\beta (wq)_{y} \eqno(34a)
$$
$$
v_{y} = k \eqno(34b)
$$
$$
w_{x} = m  \eqno(34c)
$$
where $q = \partial^{-1}_{x} k$.

\subsection{ The Myrzakulov X equation}

Similarly we can show that the M-X equation
$$
\vec S_{t} = (3k^{2} + k_{xx} + 3 \alpha^{2} w)k^{-1}\vec S_{x} \eqno(35a)
$$
$$
w_{xx} =  k_{yy} \eqno(35b)
$$
is
equivalent to the famous Kadomtsev - Petviashvili equation(KPE)
$$
(q_{t} - 6qq_{x} - q_{xxx})_{x} = 3\alpha^{2}q_{yy} \eqno(36)
$$
with $ q = k$.

\section{Conclusion}

We have discussed the differential geometrical approach to study
the some properties of soliton equations in 2+1 dimensions. Using
the geometrical formalism which was presented in[17], the
L-equivalent counterparts of the some integrable (2+1) - dimensional
spin systems are obtained. So we have shown that the Lakshmanan equivalent
soliton equations of the M-XVII, M-XI and M-X equations are the well known
modified Novikov - Veselov, Novikov - Veselov and Kadomtsev - Petviashvili
equations respectively. Finally note that the some other aspects of the some
spin systems were considered in [19-23].

\section{Appendix A: On the some integrable (2+1)-dimensional spin systems}

\subsection{ The Myrzakulov IX equation}

We note that the M-XVII equation(15) is related with
the  Myrzakulov - IX(M-IX) hierarchy
$$ iS_t+\frac{1}{2}[S,M_1S]+A_2S_x+A_1S_y = 0 \eqno(37a)$$
$$ M_2u=\frac{\alpha^{2}}{4i}tr(S[S_y,S_x]) \eqno(37b)$$
where $ \alpha,b,a  $=  consts and
$$
S= \pmatrix{
S_3 & rS^- \cr
rS^+ & -S_3
},\quad
S^{\pm}=S_{1}\pm iS_{2}, \quad S^2=I,\quad r^{2}=\pm 1
$$
$$ M_1= \alpha ^2\frac{\partial ^2}{\partial y^2}-2\alpha (b-a)\frac{\partial^2}
   {\partial x \partial y}+(a^2-2ab-b)\frac{\partial^2}{\partial x^2}; $$
$$ M_2=\alpha^2\frac{\partial^2}{\partial y^2} -\alpha(2a+1)\frac{\partial^2}
   {\partial x \partial y}+a(a+1)\frac{\partial^2}{\partial x^2},$$
$$
A_1=2i\{(2ab+a+b)u_x-(2b+1)\alpha u_y\}
$$
$$
A_2=2i\{(2ab+a+b)u_y-\alpha^{-1}(2a^2b+a^2+2ab+b)u_x\}.
$$

These set of equations is integrable and the Lax representation of
the M-IX equation(37) is given by[17]
$$ \alpha \Phi_y =\frac{1}{2}[S+(2a+1)I]\Phi_x \eqno(38a) $$
$$ \Phi_t=\frac{i}{2}[S+(2b+1)I]\Phi_{xx}+\frac{i}{2}W\Phi_x \eqno(38b) $$
with
$$ W_1=W-W_2=(2b+1)E+(2b-a+\frac{1}{2})SS_x+(2b+1)FS $$
$$ W_2=W-W_1=FI+\frac{1}{2}S_x+ES+\alpha SS_y $$
$$ E = -\frac{i}{2\alpha} u_x,\,\,\,  F = \frac{i}{2}(\frac{(2a+1)u_{x}}{\alpha} -
2u_{y}) $$
It is well known that the M-IX equation(37) is equivalent to the following
Zakharov equation[13]
$$
iq_{t}+M_{1}q+vq=0, \eqno(39a)
$$
$$
ip_{t}-M_{1}p-vp=0, \eqno(39b)
$$
$$
M_{2}v = -2M_{1}(pq), \eqno(39c)
$$
where $v=i(c_{11}-c_{22}),\quad p= E\bar q$.
It is interest note that the M-IX equation(37) admits the some
integrable reductions. Let us now present these particular integrable
cases.

\subsection{ The Myrzakulov VIII equation}

Let $b=0$. Then equations(37) take the form
$$ iS_t=\frac{1}{2}[S_{\xi\xi},S]+iwS_{\xi}   \eqno (40a) $$
$$ w_{\eta}=\frac{1}{4i}tr(S[S_{\xi},S_{\eta}])  \eqno (40b) $$
where
$$
\xi = x+\frac{a+1}{\alpha}y,\,\,\,\,\,
\eta = -x -\frac{a}{\alpha}y,\,\,\,\,\, w=u_{\xi},
$$
which is the M-VIII equation[17]. The  equivalent counterpart
of the M-VIII equation(40) we obtain from(39) as $b=0$
$$
iq_{t}+q_{\xi \xi}+vq=0, \eqno(41a)
$$
$$
v_{\eta} = -2r^{2}(\bar q q)_{\xi}, \eqno(41b)
$$
which is the other Zakharov equation[13].

\subsection{ The Ishimori equation}

Now consider the case: $ a=b=-\frac{1}{2} $. In this case
equations(37) reduces to the well known Ishimori equation
$$ iS_t+\frac{1}{2}[S,(\frac{1}{4}S_{xx}+\alpha^{2}S_{yy})]+
iu_{y}S_x+iu_{x}S_y = 0 \eqno(42a)$$
$$ \alpha^{2}u_{yy} - \frac{1}{4}u_{xx}=
\frac{\alpha^{2}}{4i}tr(S[S_y,S_x]) \eqno(42b)$$
The  equivalent counterpart of the equation(42) is the Davey-Stewartson
equation
$$ iq_t+\frac{1}{4}q_{xx}+\alpha^{2}q_{yy}+
vq = 0 \eqno(43a)$$
$$ \alpha^{2}v_{yy} - \frac{1}{4}v_{xx}=-2\{\alpha^{2}(pq)_{yy}
 +\frac{1}{4}(pq)_{xx}\} \eqno(43b)$$
that follows from the ZE(39). Note that equations (42) and (43) are
gauge equivalent each to other[11].

\subsection{ The Myrzakulov XVIII equation}

Consider the reduction: $a=-\frac{1}{2}$.
Then (37) reduces to the M-XVIII
equation[17]
$$ iS_t+\frac{1}{2}[S,(\frac{1}{4}S_{xx}-\alpha(2b+1)S_{xy}+\alpha^{2}S_{yy})]+
A_{20}S_x+A_{10}S_y = 0 \eqno(44a)$$
$$ \alpha^{2}u_{yy} - \frac{1}{4}u_{xx}=\frac{\alpha^{2}}{4i}tr(S[S_y,S_x]) \eqno(44b)$$
where $A_{j0}=A_{j}$ as $a=-\frac{1}{2}$. The corresponding gauge equivalent equation
obtain from(39) and looks like
$$ iq_t+\frac{1}{4}q_{xx}-\alpha(2b+1)q_{xy}+\alpha^{2}q_{yy}+
vq = 0 \eqno(45a)$$
$$ \alpha^{2}v_{yy} - \frac{1}{4}v_{xx}=-2\{\alpha^{2}(pq)_{yy}-
\alpha (2b+1)(pq)_{xy} +\frac{1}{4}(pq)_{xx}\} \eqno(45b)$$
Note that the Lax representations of equations(40), (42) and (44)
we can get from (38) as $b=0, a=b=-\frac{1}{2}$ and $a=-\frac{1}{2}$
respectively.

\section{Appendix B: Differential
geometry of surfaces and the M-XXII equation}

Consider the Myrzakulov (M-XXII) equation[17]
$$
-iS_t=\frac{1}{2}([S,S_y]+2iuS)_x+\frac{i}{2}V_1S_x-2ib^2 S_y \eqno(46a)
$$
$$
u_x=-\vec S(\vec S_x\wedge \vec S_y), \quad
V_{1x}=\frac{1}{4b^2}(\vec S^2_x)_y, \eqno(46b)
$$
where $\vec S = (S_{1}, S_{2}, S_{3})$ is a spin vector,
$\vec S^{2} = E = \pm 1$.
These equations are integrable. The corresponding Lax representation is
given by[17]
$$
\Phi_{x} = \{ - i(\lambda^{2} - b^{2})S + \frac{\lambda - b}{2b}SS_{x}\}\Phi  \eqno(47a)
$$
$$
\Phi_{t} = 2\lambda^{2}\Phi_{y} + \{(\lambda^{2} - b^{2})(2A + B) +
(\lambda - b)C\}\Phi  \eqno(47b)
$$
with
$$
A = \frac{1}{4}([S,S_{y}] +2iuS) + \frac{i}{4}V_{1}S,\quad
B = \frac{i}{2}V_{1}S, \quad
C = - \frac{V_{1}}{4b^{2}}SS_{x} +\frac{i}{2b}\{ S_{xy} - [S_{x}, A]\}.
$$
Here
spin
matrix has the form
$$
S= \pmatrix{
S_3 & rS^- \cr
rS^+ & -S_3
}, \quad S^2=I,\quad r^2=\pm1,\quad S^{\pm } = S_{1} \pm i S_{2}.
$$

Now find the Lakshmanan equivalent counterpart of the M-XXII equation (46)
for the case $E = +1$ (for the case $E = -1$, see, e.g., [17]).
To this end  we can use the two geometrical approaches(D- and C-approaches).
Let us use the C-approach, i.e., the surface
approach. Consider the motion of surface in
the 3-dimensional space which generated by a position vector
$\vec r(x,y,t) = \vec r(x^{1}, x^{2}, t)$. According to the C-approach,
$x$ and $y$ are local coordinates on the surface. The
first and second fundamental forms in the usual notation are given by
$$ I=d\vec r d\vec r=Edx^2+2Fdxdy+Gdy^2, \quad
II=-d\vec r d\vec n=Ldx^2+2Mdxdy+Ndy^2 \eqno (48) $$
where $$ E=\vec r_x\vec r_x=g_{11},\quad  F=\vec r_x\vec r_y=g_{12},\quad
G=\vec r_{y^2}=g_{22},   $$
$$ L=\vec n\vec r_{xx}=b_{11},\quad M=\vec n\vec r_{xy}=b_{12},\quad
N=\vec n\vec r_{yy}=b_{22},\quad \vec n=\frac{(\vec r_x \wedge \vec r_y)}.
{|\vec r_x \wedge \vec r_y|}.  $$

In this case, the set of equations of the C-approach[17],
becomes
$$ \vec r_{t} = W_{1}\vec r_x + W_{2} \vec r_y + W_{3} \vec n \eqno (49a) $$
$$ \vec r_{xx}=\Gamma^1_{11} \vec r_x + \Gamma^2_{11} \vec r_y +L \vec n \eqno (49b) $$
$$ \vec r_{xy}=\Gamma^1_{12} \vec r_x + \Gamma^2_{12} \vec r_y +M \vec n \eqno (49c)$$
$$ \vec r_{yy}=\Gamma^1_{22} \vec r_x + \Gamma^2_{22} \vec r_y + N \vec n \eqno (49d)$$
$$ \vec n_x=p_1 \vec r_x+p_2 \vec r_y \eqno (49e) $$
$$ \vec n_y=q_1 \vec r_x+q_2 \vec r_y \eqno (49f) $$
where $W_{j}$ are some functions, $ \Gamma^k_{ij} $ are the Christoffel symbols of the second kind defined by
the metric $ g_{ij} $ and $ g^{ij}=(g_{ij})^{-1} $ as
$$ \Gamma^k_{ij}=\frac{1}{2} g^{kl}(\frac{\partial g_{lj}} {\partial x^i}+
   \frac {\partial g_{il}}{ \partial x^j}-\frac{\partial g_{ij}}
   {\partial x^l}) \eqno (50) $$
The coefficients $ p_i, q_i $ are given by
$$ p_i=-b_{1j}g^{ji}, \,\,\, q_i=-b_{2j}g^{ji}. \eqno (51) $$
The compatibility conditions $ \vec r_{xxy}=\vec r_{xyx} $ and
$ \vec r_{yyx}=\vec r_{xyy} $ yield the following Mainardi-Peterson-Codazzi
equations (MPCE)
$$ R^l_{ijk} = b_{ij}b^l_{k}-b_{ik}b^l_{j},\quad
\frac{\partial b_{ij}}{\partial x^k}-\frac{\partial b_{ik}}{\partial
    x^j}=\Gamma^s_{ik}b_{is}-\Gamma^s_{ij}b_{ks}
\eqno (52) $$
where $ b^j_i=g^{jl}b_{il} $ and the curvature tenzor has the form
$$ R^l_{ijk} = \frac{\partial \Gamma^l_{ij}}{\partial x^k}-\frac{\partial \Gamma^l_{ik}}
{\partial x^j}+\Gamma^s_{ij} \Gamma^l_{ks}-\Gamma^s_{ik} \Gamma^l_{js}
\eqno (53) $$

Let
$ Z = ( r_{x},  r_{y},  n)^{t}$ . Then
$$ Z_{x} = A  Z,\quad Z_{y} = B  Z \eqno (54) $$
where
$$
A =
\left ( \begin{array}{ccc}
\Gamma^{1}_{11} & \Gamma^{2}_{11} & L \\
\Gamma^{1}_{12} & \Gamma^{2}_{12} & M \\
p_{1}           & p_{2}           & 0
\end{array} \right), \quad
B =
\left ( \begin{array}{ccc}
\Gamma^{1}_{12} & \Gamma^{2}_{12} & M \\
\Gamma^{1}_{22} & \Gamma^{2}_{22} & N \\
q_{1}           & q_{2}           & 0
\end{array} \right). \eqno(55)
$$
Hence we get the new form of the MPCE(52)
$$ A_y - B_x + [A, B] = 0 \eqno (56) $$

Let us introduce the orthogonal trihedral[17]
$$ \vec e_{1} = \frac{\vec r_x}{\surd E}, \,\,\,
\vec e_{2} = \vec n, \,\,\, \vec e_{3} = \vec e_{1} \wedge
\vec e_{2}    \eqno(57) $$

Let $ \vec r_x^2 = E = \pm 1 $ and $ F = 0$.  Then the vectors $
\vec e_{j}$   satisfy the following Serret-Frenet-Myrzakulov equation(SFME)
$$
\left ( \begin{array}{ccc}
\vec e_{1} \\
\vec e_{2} \\
\vec e_{3}
\end{array} \right)_{x} = C
\left ( \begin{array}{ccc}
\vec e_{1} \\
\vec e_{2} \\
\vec e_{3}
\end{array} \right),\quad
\left ( \begin{array}{ccc}
\vec e_{1} \\
\vec e_{2} \\
\vec e_{3}
\end{array} \right)_{y} = D
\left ( \begin{array}{ccc}
\vec e_{1} \\
\vec e_{2} \\
\vec e_{3}
\end{array} \right) \eqno(58)
$$
and the time equation
$$
\left ( \begin{array}{ccc}
\vec e_{1} \\
\vec e_{2} \\
\vec e_{3}
\end{array} \right)_{t} = G
\left ( \begin{array}{ccc}
\vec e_{1} \\
\vec e_{2} \\
\vec e_{3}
\end{array} \right) \eqno(59)
$$
with
$$
C =
\left ( \begin{array}{ccc}
0   & k     & 0 \\
-Ek & 0     & \tau  \\
0   & -\tau & 0
\end{array} \right) ,\quad
D =
\left ( \begin{array}{ccc}
0       & m_{3}  & -m_{2} \\
-Em_{3} & 0      & m_{1} \\
Em_{2}  & -m_{1} & 0
\end{array} \right),\quad
G =
\left ( \begin{array}{ccc}
0       & \omega_{3}  & -\omega_{2} \\
-E\omega_{3} & 0      & \omega_{1} \\
E\omega_{2}  & -\omega_{1} & 0
\end{array} \right),
$$
where
$$k =\frac{L}{2}, \quad
\tau = MG^{-1/2}.     \eqno(60)
$$
Hence we have
$$C_y - D_x + [C, D] = 0,\quad
C_t - G_x + [C, G] = 0,\quad
D_t - G_y + [D, G] = 0.  \eqno(61)
$$
Now let $\vec S = \vec e_{1}$.
Let us now introduce the new function $q$ by
$$
q = \frac{k}{2b}\exp{i[\frac{1}{8}\partial^{-1}_{x}(k^{2}b^{-2} - 4\tau)
- 2b^{2}x]}  \eqno(62)
$$
Then the function $q$ satisfies the following equations[17]
$$
iq_t + q_{yx} + \frac{i}{2}[(V_1q)_x - V_2 q - qpq_y] = 0 \eqno(63a)
$$
$$
ip_t - p_{yx} + \frac{i}{2}[(V_1p)_x + V_{2}p - qpp_y] = 0 \eqno(63b) $$
$$ V_{1x}=(pq)_y, \quad  V_{2x}=p_{yx}q-pq_{yx} \eqno(63c)
$$
where $p = E\bar q$. This  set of equations is the
L-equivalent counterpart
of the M-XXII equation(46). As it seems to us, equations (63) are new integrable
equation.
We will call (63) the M-XXII$_{q}$ equation.
Now let us consider the following transformation
$$
q^{\prime} = q\exp(-\frac{i}{2}\partial^{-1}_{x}|q|^{2})   \eqno(64)
$$

Then the new variable $q^{\prime}$ satisfies the Strachan equation[18]
$$ iq^{\prime}_{t} + q^{\prime}_{xy} + i(Vq^{\prime})_x = 0, \quad
 V_x = E(|q^{\prime}|^2)_y. \eqno (65) $$

We see that  the M-XXII$_{q}$ equation(7) and the Strachan equation(22) is
gauge equivalent to each other.

\end{document}